\documentclass[conference,final]{IEEEtran}

\usepackage{latexsym}
\usepackage{amsfonts}
\usepackage{amssymb}
\usepackage{graphicx}
\usepackage{amscd}
\usepackage{amsmath}
\usepackage{amsfonts}
\usepackage{amssymb}

\newcommand{\RR}{\mathbb{R}}
\newcommand{\p}{\rm {I\kern-1pt P}}

\newtheorem{theorem}{Theorem}[section]
\newcommand{\EQ}{\begin{equation}\begin{array}{lllllllll}}
\newcommand{\EE}{\end{array}\end{equation}}
\newcommand{\MT}{\left[ \begin{array}{ccccccccc}}
\newcommand{\EM}{\end{array}\right]}

\newcommand{\LM}{\left[\begin{array}{ccccccccc}}
\newcommand{\RM}{\end{array}\right]}
\newcommand{\LA}{\left\{ \begin{array}{ccccccccc}}
\newcommand{\RA}{\end{array}\right.}
\newcommand{\RAA}{\end{array}\right\}}

\newcommand{\pfrac}[2]{\frac{\partial#1}{\partial#2}}

\newcommand{\bc}{\mathbf{c}}
\newcommand{\bC}{\mathbf{C}}
\newcommand{\bk}{\mathbf{k}}
\newcommand{\cH}{\mathcal{H}}

\newcommand{\cX}{\mathcal{X}}
\newcommand{\Tx}{C_{\mathbf{x}}}
\newcommand{\tr}{\top\!}
\newcommand{\scal}[2]{\left\langle{#1},{#2}\right\rangle}
\newcommand{\nor}[1]{\left\|{#1}\right\|}

\def\[{\begin{equation}}
\def\]{\end{equation}}

\def\Box{\vrule width1ex height1ex}

\def\[{\begin{equation}}
\def\]{\end{equation}}

\def\Box{\vrule width1ex height1ex}

\newcommand{\eq}{\begin{equation}\begin{array}{lclllllllllllllll}}
\newcommand{\ee}{\end{array}\end{equation}}
\newcommand{\bmt}{\left[ \begin{array}{ccccccccc}}
\newcommand{\emt}{\end{array}\right]}
\newcommand{\bea}{\begin{eqnarray}}
\newcommand{\eea}{\end{eqnarray}}
\newcommand{\bean}{\begin{eqnarray*}}
\newcommand{\eean}{\end{eqnarray*}}

\begin{document}

\title{\bf\LARGE Balanced Reduction of Nonlinear Control Systems in Reproducing Kernel Hilbert Space}
\author{\IEEEauthorblockN{Jake Bouvrie}
\IEEEauthorblockA{Department of Mathematics\\
Duke University\\
Durham, NC 27708 USA\\
Email: jvb@math.duke.edu}
\and
\IEEEauthorblockN{Boumediene Hamzi}
\IEEEauthorblockA{Department of Mathematics\\
Duke University\\
Durham, NC 27708 USA\\
Email: hamzi@math.duke.edu}}

\maketitle
\begin{abstract}
We introduce a novel data-driven order reduction method for nonlinear control systems, drawing
on recent progress in machine learning and statistical dimensionality reduction.
The method rests on the assumption that the nonlinear system behaves linearly when
lifted into a high (or infinite) dimensional feature space where balanced truncation
may be carried out implicitly. This leads to a nonlinear reduction map which can be combined
with a representation of the system belonging to a
reproducing kernel Hilbert space to give a closed, reduced order dynamical system
which captures the essential input-output characteristics of the original model.
Empirical simulations illustrating the approach are also provided.
\end{abstract}

\section{Introduction}
Model reduction of controlled dynamical systems has been a long standing,
and as yet, unsettled challenge in control theory. The benefits are clear:
a low dimensional approximation of a high dimensional system can be manipulated
with a simpler controller, and can be simulated at lower computational cost.
A complex, high dimensional system may even be replaced by a simpler model
all together leading to significant cost savings, as in circuit design, while
the ``important variables'' of a system might shed light on underlying physical
or biological processes.
Reduction of linear dynamical systems has been treated
with some success to date. As we describe in more detail below, model reduction in the linear
case proceeds by reducing the dimension of the system with an eye towards preserving
its essential input-output behavior, a notion directly related to ``balancing'' observability and
controllability of the system. The nonlinear picture, however, is considerably more
involved.

In this paper we propose a scheme for balanced model-order reduction of general, nonlinear
control systems. A key, and to our knowledge, novel point of departure from the literature
on nonlinear model reduction is that our approach marries approximation and dimensionality reduction methods known to the machine learning and statistics communities with existing ideas in linear and nonlinear control. In particular, we apply a method similar to kernel PCA as well as function learning in Reproducing Kernel Hilbert Spaces (RKHS) to the problem of balanced model reduction. Working in RKHS
provides a convenient, general functional-analytical framework for theoretical understanding
as well as a ready source of existing results and error estimates. The approach presented here is also strongly empirical, in that observability and controllability, and in some cases the dynamics of the nonlinear system are estimated from simulated or measured trajectories. This emphasis on the empirical makes our approach broadly applicable, as the method can be applied without having to tailor anything to the particular form of the dynamics.

The approach we propose begins by constructing empirical estimates of the observability and controllability
Gramians in a high (or possibly infinite) dimensional feature space. The Gramians are simultaneously
diagonalized in order to identify directions which, in the feature space, are both the most
observable and the most controllable. The assumption that a nonlinear system behaves linearly when lifted to a feature space is far more reasonable than assuming linearity in the original space, and then carrying out the linear theory hoping for the best. Working in the high dimensional feature space allows one to
perform linear operations on a representation of the system's state and output which can capture strong nonlinearities. Therefore a system which is not model reducible using existing methods, may become reducible when mapped into such a nonlinear feature space. This situation closely parallels the problem of linear separability in data classification: A dataset which is not linearly separable might be easily
separated when mapped into a nonlinear feature space. The decision boundary is linear in this feature space, but is nonlinear in the original data space.

Nonlinear reduction of the state space already opens the door to the design of simpler controllers,
but is only half of the picture. One would also like to be able to write a closed, reduced dynamical
system whose input-output behavior closely captures that of the original system. This problem is the
focus of the second half of our paper, where we again exploit helpful properties of RKHS in order to provide such a closed system.

The paper is organized as follows. In the next section we provide the relevant background for
model reduction and balancing. We then adapt and extend balancing techniques described in the background to the current RKHS setting in Section~\ref{sec:empirical_gramians}. Section~\ref{sec:closed_sys} then
proposes a method for determining a closed, reduced nonlinear control system in light of the reduction
map described in~Section~\ref{sec:empirical_gramians}.
Finally,
 Section~\ref{sec:expts} provides experiments illustrating an application of the proposed methods to
a specific nonlinear system.

\section{Background}
Several approaches have been proposed for the reduction of
linear control systems in view of control, but few
exist for finite or infinite-dimensional controlled
nonlinear dynamical systems. For linear systems the pioneering ``Input-
Output balancing'' approach proposed by B.C. Moore observes that the important states are the
ones that are both easy to reach and that generate a lot of energy at the output. If a large amount of energy is required to reach a certain state but the same state yields a small output energy, the state is unimportant for the input-output behavior of the system.
The goal is then to find the states that are \emph{both} the most controllable and the most observable. One way to determine such states is to find a change of coordinates where the controllability and observability Gramians (which can be viewed as a measure of the controllability and the observability of the system) are equal and diagonal. States that are
difficult to reach and that don't significantly affect the output are then ignored or truncated. A system expressed in the coordinates where each state is equally controllable and observable is called its \emph{balanced realization}.

A proposal for generalizing this approach to nonlinear control systems was advanced by J. Scherpen~\cite{scherpen_thesis}, where suitably defined controllability and observability energy functions reduce to Gramians in the linear case. In general, to find the balanced realization of a system one needs to solve a set of Hamilton-Jacobi and Lyapunov equations (as we will discuss below). Moore~\cite{moore} proposed an alternative, data-based approach for balancing in the linear case. This method uses samples of the impulse response of a linear system to construct empirical controllability and observability Gramians which are then balanced and truncated using Principal Components Analysis (PCA, or POD). This data-driven strategy was then extended to nonlinear control systems with a stable linear approximation by Lall et al.~\cite{lall}, by effectively applying Moore's method to a nonlinear system by way of the Galerkin projection. Despite the fact that the balancing theory underpinning their approach assumes a linear system, Lall and colleagues were able to effectively reduce some nonlinear systems.

Phillips~\cite{Phillips} et al. has also studied reduction of nonlinear circuit models
in the case of linear but unbalanced coordinate transformations and found
that approximation using a polynomial RKHS
could afford computational advantages. Gray and Verriest mention
in~\cite{gray} that studying algebraically defined Gramian operators
in RKHS may provide advantageous approximation properties, though the idea
is not further explored. Finally, Coifman et al.~\cite{Mauro08} discuss
 reduction of an uncontrolled stochastic Langevin system. There, eigenfunctions of a
combinatorial Laplacian, built from samples of trajectories, provide a
set of reduction coordinates but does not provide a reduced system. This method is
related to kernel principal components (KPCA) using a Gaussian kernel, however
reduction in this study is carried out on a simplified linear
system outside the context of control.

In the following section we review balancing
of linear and nonlinear systems as introduced in~\cite{moore} and~\cite{scherpen_thesis}.

\subsection{Balancing of Linear Systems}
Consider a linear control system \[\label{linsys}\begin{array}{rcl}\dot{x}&=&Fx+Gu,\\y&=&Hx, \end{array}, \]
where $(F,G)$ is controllable, $(F,H)$ is observable and $F$ is Hurwitz. We define the controllability and the observability Gramians as, respectively,
\[ \label{gramians}\begin{array}{rcl}
W_c&=&\int_0^{\infty}e^{Ft}GG^{\tr}e^{F^{\tr}t}\, dt,\\
W_o&=&\int_0^{\infty}e^{F^{\tr}t}H^{\tr}He^{Ft}\, dt.
\end{array}\nonumber\]
These two matrices can be viewed as a measure of the controllability and the observability of the system~\cite{moore}. For instance, consider the past energy~\cite{scherpen_thesis}, $L_c(x_0)$, defined
as the minimal energy required to reach $x_0$ from $0$ in infinite time
\[\label{L_c}
L_c(x_0)=\inf_{\substack{
u \in { L}_2(-\infty,0),\\ x(-\infty)=0, x(0)=x_0}}
\frac{1}{2}\int_{-\infty}^0||u(t)||^2\,dt,
\]
 and  the future energy~\cite{scherpen_thesis}, $L_o(x_0)$, defined as the output energy generated
by releasing the system from its initial state $x(t_0)=x_0$, and zero input $u(t)=0$ for $t\ge0$, i.e.
 \[\label{L_o}
 L_o(x_0)=\frac{1}{2}\int_{0}^{\infty}||y(t)||^2\,dt,
 \]
for $x(t_0)=x_0$ and $u(t)=0, t\ge0$.
In the linear case, it can be shown that $L_c(x_0)=\tfrac{1}{2}x_0^{\tr}W_c^{-1}x_0,$ and $L_o(x_0)=\tfrac{1}{2}x_0^{\tr}W_o x_0.$ The columns of $W_c$ span the controllable subspace while the nullspace of $W_o$ coincides with the unobservable subspace. As such, $W_c$ and $W_o$ (or their estimates) are the key ingredients in many model reduction techniques.
It is also well known that $W_c$ and $W_o$ satisfy the Lyapunov equations~\cite{moore}
\[\label{gramians_lyap}
\begin{array}{rcl}
FW_c+W_cF^{\tr}=-GG^{\tr},\\
F^{\tr}W_o+W_oF=-H^{\tr}H .
\end{array}
\nonumber\]
Several methods have been developed to solve these equations directly~\cite{laub,li}.

The idea behind balancing is to find a representation where the system's observable and controllable subspaces are aligned so that reduction, if possible, consists of eliminating uncontrollable states which are also the least observable. More formally, we would like to find a new coordinate system such that
\[W_c=W_o=\Sigma=\mbox{diag}\{\sigma_1,\cdots,\sigma_n\},\nonumber\]
where $\sigma_1 \ge \sigma_2 \ge \cdots \ge \sigma_n \ge 0$.
If $(F,G)$ is controllable and $(F,H)$ is observable, then there exists a transformation such that the  state space expressed in the transformed coordinates $(TFT^{-1},TG,HT^{-1})$ is
balanced and $TW_cT^{\tr}=T^{-{\tr}}W_oT^{-1}=\Sigma$. Typically one looks for a gap in the singular values $\{\sigma_i\}$ for guidance as to where truncation should occur. If we see that there is a  $k$ such that $\sigma_k \gg\sigma_{k+1}$, then the states most responsible for governing the input-output relationship of the system are $(x_1,\cdots,x_k)$ while $(x_{k+1},\ldots,x_n)$ are assumed to make negligible contributions.

Although several methods exist for computing $T$~\cite{laub,li}, the general idea is to compute the Cholesky decomposition of $W_o$ so that $W_o=Z Z^{\tr}$, and form the SVD $U \Sigma^2 U^{\tr}$ of  $Z^{\tr} W_c Z$. Then $T$ is given by $T=\Sigma^{\frac{1}{2}}U^{\tr}Z^{-1}$.
We also note that the problem of finding the coordinate change $T$ can be seen as an optimization problem~\cite{antoulas} of the form
\[\min_{T} \mbox{trace}[TW_cT^{\ast}+T^{-\ast}W_oT^{-1}]. \nonumber\]

\subsection{Balancing of Nonlinear Systems}
In the nonlinear case, the energy functions $L_c$ and $L_o$ in (\ref{L_c}) and (\ref{L_o}) are obtained by solving both a Lyapunov and a Hamilton-Jacobi equation. Here we follow the development of Scherpen~\cite{scherpen_thesis}. Consider the nonlinear system
\[\label{sigma}
\left\{\begin{array}{rcl}\dot{x}&=&f(x)+\sum_{i=1}^mg_i(x)u_i,\\ y &=& h(x), \end{array}\right.
\]
with $x \in \RR^n$, $u \in \RR^m$, $y\in \RR^p$, $f(0)=0$, $g_i(0)=0$ for $1 \le i \le m$, and $h(0)=0$.  Moreover, assume the following Hypothesis.\\
{\it Hypothesis H:} The linearization of~\eqref{sigma} around the origin is controllable, observable and $F=\frac{\partial f}{\partial x}|_{x=0}$ is asymptotically stable.
\begin{theorem}\label{thm:scherp1}\cite{scherpen_thesis}
 If the origin is an asymptotically stable equilibrium of $f(x)$ on a neighborhood $W$ of the origin, then for all $x \in W$, $L_o(x)$ is the unique smooth solution of
\[\label{Lo_hjb} \frac{\partial L_o}{\partial x}(x)f(x)+\frac{1}{2}h^{\tr}(x)h(x)=0,\quad L_o(0)=0 \]
under the assumption that (\ref{Lo_hjb}) has a smooth solution on $W$. Furthermore for all $x \in W$, $L_c(x)$ is the unique smooth solution of
\[\label{Lc_hjb} \frac{\partial L_c}{\partial x}(x)f(x)+\frac{1}{2} \frac{\partial L_c}{\partial x}(x)g(x)g^{\tr}(x)  \frac{\partial^{\tr} L_c}{\partial x}(x)=0, \quad L_c(0)=0\]
under the assumption that (\ref{Lc_hjb}) has a smooth solution $\bar{L}_c$ on $W$ and that the origin is an asymptotically stable equilibrium of $-(f(x)+g(x)g^{\tr}(x) \frac{\partial \bar{L}_c}{\partial x}(x))$ on $W$.
\end{theorem}
With the controllability and the observability functions on hand, the input-normal/output-diagonal realization of system~\eqref{sigma} can be computed by way of a coordinate transformation. More precisely,
\begin{theorem}\label{theorem_scherpen}\cite{scherpen_thesis}
Consider system~\eqref{sigma} under Hypothesis H and the assumptions in Theorem~\ref{thm:scherp1}. Then, there exists a neighborhood $W$ of the origin and coordinate transformation $x=\varphi(z)$ on $W$ converting  the energy functions  into the form
\[L_c(\varphi(z))=\frac{1}{2}z^{\tr}z, \nonumber\]
\[L_o(\varphi(z))=\frac{1}{2}\sum_{i=1}^nz_i^2\sigma_i(z_i)^2, \nonumber\]
where $\sigma_1(x) \ge \sigma_2(x) \ge \cdots \ge \sigma_n(x)$. The functions $\sigma_i(\cdot)$ are called {\em Hankel singular value functions}.
\end{theorem}
Analogous to the linear case, the system's states can be sorted in order of importance by sorting
the singular value functions, and reduction proceeds by removing the least important states.

In the above framework for balancing of nonlinear systems, one needs to solve (or numerically evaluate) the PDEs (\ref{Lo_hjb}), (\ref{Lc_hjb}) and compute the coordinate change $x=\varphi(z)$, however there are no systematic methods or tools for solving these problems. Various approximate solutions based on Taylor series expansions have been proposed~\cite{krener1,krener2,fujimoto}. Newman~\cite{newman} introduces a statistical approximation based on exciting the system with white Gaussian noise and then computing the balancing transformation using an algorithm from differential topology. As mentioned earlier, an essentially linear empirical approach was proposed in~\cite{lall}. In this paper, we combine aspects of both data-driven approaches and analytic approaches by carrying out balancing in a suitable RKHS.

\section{Empirical Balancing of Nonlinear Systems in RKHS}\label{sec:empirical_gramians}
We consider a general nonlinear system of the form
\[\label{eqn:nlsys}
\left\{\begin{array}{rcl} \dot{x}&=&f(x,u)\\ y&=&h(x) \end{array}\right.
\]
with $x \in \RR^n$, $u \in \RR^m$, $y \in \RR^p$, $f(0,0)=0$, and $h(0)=0$. Let ${\cal R}(x_0)=\{x' \in \RR^n: \exists\, u \in L_{\infty}(\RR,\RR^m) \;\;\mbox{and}\;\; \exists\, T \in [0,\infty) \;\; \mbox{such that}\; \;  x(0)=x_0\;\; \mbox{and}\;\; x(T)=x' \}$ be the reachable set from the initial condition $x(0)=x_0$.
We assume that the system is zero-state observable, and that the linearization of (\ref{eqn:nlsys}) around the origin is controllable. We also assume that the origin of $\dot{x}=f(x,0)$ is asymptotically stable.

We treat the problem of estimating the observability and controllability Gramians as one of estimating an integral operator from data in a reproducing kernel Hilbert space (RKHS)~\cite{AronRKHS}. \emph{Our approach hinges on the key modeling assumption that the nonlinear dynamical system is linear in an appropriate high (or possibly infinite) dimensional lifted feature space}. Covariance operators in this feature space and their empirical estimates are the objects of primary importance and contain the information needed to perform model reduction. In particular, the (linear) observability and controllability Gramians  are estimated and diagonalized in the feature space, but capture nonlinearities in the original state space. The reduction approach we propose adapts ideas from kernel PCA (KPCA)~\cite{KPCA:98} and is driven by a set of simulated or sampled system trajectories, extending and generalizing the work of Moore~\cite{moore} and Lall et al.~\cite{lall}.


\subsection{Definitions}\label{sec:defins}
In the development below we lift state vectors of the system into a reproducing kernel Hilbert space~\cite{AronRKHS}, $\cH$, endowed with a symmetric positive definite kernel function $K:\cX\times\cX\to\RR$ which we assume here to be continuous and bounded by $\kappa = \sup_{x\in\cX}\sqrt{K(x,x)} < \infty$. In particular, we make use of the following important properties: For all $f\in\cH$, $f(x) = \scal{f}{K_x}_{\cH}$, where $K_x:=K(x,\cdot)$. This is the reproducing property. Second, to any RKHS we can associate a feature map $\Phi:\cX\to\mathcal{F}$ satisfying $\scal{\Phi(x)}{\Phi(x')}_{\cH}=K(x,x')$. For example, we can take $\Phi(x):=K_x$ in which case $\mathcal{F}=\cH$ -- the ``feature space'' is the RKHS. We will further assume that $\cH$ is always separable. 

\subsection{Empirical Gramians}
Following~\cite{moore}, we estimate the controllability Gramian by exciting each coordinate
of the input with impulses while setting $x_0 = 0$. One can also further excite using rotations of impulses as suggested in~\cite{lall}, however for simplicity we consider only the original signals proposed in~\cite{moore}. Let $u^i(t) = \delta(t)e_i$ be the $i$-th excitation
signal, and let $x^i(t)$ be the corresponding response of the system. Form the matrix $X(t) = \bigl[x^1(t) ~\cdots~ x^m(t)\bigr] \in \RR^{n\times m}$, so that $X(t)$ is seen as a data matrix with
column observations given by the respective responses $x^i(t)$. Then $W_c\in\RR^{n\times n}$ is given by
\[
W_c = \frac{1}{m}\int_0^{\infty}X(t)X(t)^{\tr} .
\nonumber\]
We can approximate this integral by sampling the matrix function $X(t)$ within a finite time interval $[0,T]$ assuming the regular partition $\{t_i\}_{i=1}^N, t_i = (T/N)i$. This leads to the empirical controllability Gramian
\[
\widehat{W}_c = \frac{T}{mN}\sum_{i=1}^N X(t_i)X(t_i)^{\tr} .
\nonumber\]

As described in~\cite{moore}, the observability Gramian is estimated by
fixing $u(t) = 0$, setting $x_0 = e_i$ for $i=1,\ldots,n$, and measuring the corresponding system output responses $y^i(t)$. As before, assemble the responses into a matrix $Y(t) = [y^1(t) ~\cdots~ y^n(t)]\in \RR^{p\times n}$. The observability Gramian $W_o\in\RR^{n\times n}$ and its empirical
counterpart $\widehat{W}_o$ are given by
\[
W_o = \frac{1}{p}\int_0^{\infty}Y(t)^{\tr}Y(t)\, , \quad \widehat{W}_o = \frac{T}{pN}\sum_{i=1}^N \widetilde{Y}(t_i)\widetilde{Y}(t_i)^{\tr}
\nonumber\]
where $\widetilde{Y}(t) = Y(t)^{\tr}$.
The matrix $\widetilde{Y}(t_i)\in\RR^{n\times p}$ can be thought of as a data matrix with column observations
\begin{equation}\label{eqn:obs_data}
d_j(t_i) = \bigl(y_j^1(t_i), \ldots, y_j^n(t_i)\bigr)^{\!\tr} \in\RR^n,\quad j=1,\ldots,p,
\end{equation}
so that $d_j(t_i)$ corresponds to the response at time $t_i$ of the single output coordinate $j$ to each of the (separate) initial conditions $x_0=e_k, k=1,\ldots,n$. This convention will lead to greater clarity in the steps that follows.

\subsection{Kernel PCA}
Kernel PCA~\cite{KPCA:98} generalizes linear PCA by carrying out PCA in a high dimensional feature space defined by
a feature map $\Phi:\RR^{n}\to\mathcal{F}$. Taking the feature map $\Phi(x)=K_x$ and given the set of data $\mathbf{x}:=\{x_i\}_{i=1}^N \in\RR^{n}$, we can consider PCA in the feature space by simply working with the covariance of the mapped vectors,
\begin{equation}\label{eqn:kpca_cov}
\Tx = \frac{1}{N}\sum_{i=1}^N\Phi(x_i)\otimes \Phi(x_i),
\end{equation}
where $\Phi(x_i)\otimes \Phi(x_i) = \scal{\Phi(x_i)}{\cdot}\Phi(x_i)$ denotes the tensor product between two vectors in $\cH$.
We will assume the data are centered in the feature space so that $\sum_i \Phi(x_i) = 0$. If not, data may be centered according to the prescription in~\cite{KPCA:98}.
The principal subspaces are computed by diagonalizing $\Tx$, however as is shown in~\cite{KPCA:98}, one can equivalently form the matrix $K\in\RR^{N\times N}$ of kernel
products $(K)_{ij} = K(x_i, x_j)$ for $i,j=1,\ldots,N$, and solve the eigenproblem $K\boldsymbol{\alpha}=N\lambda\boldsymbol{\alpha}$. If $\Tx v_i = \lambda_i v_i$, then we have that $v_i = \Psi\boldsymbol{\alpha}_i$ where
$\Psi:=\bigl(\Phi(x_1)~\cdots~\Phi(x_N)\bigr)$, and the non-zero eigenvalues of $K$ and $\Tx$ coincide.
The eigenvectors $\boldsymbol{\alpha}_i$ of $K$ are then normalized so that the eigenvectors $v_i$ of $\Tx$ have unit norm in the feature space, leading to the condition $\nor{\boldsymbol{\alpha}_i}^2=\lambda_i^{-1}$. Assuming this normalization convention, sort the eigenvectors according to the magnitudes of the corresponding eigenvalues in descending order, and form the matrix $A_q = \bigl[\boldsymbol{\alpha}_1 ~\cdots~ \boldsymbol{\alpha}_q\bigr], 1\leq q\leq \min(n,N)$. Similarly, form the matrix $V_q = \bigl[v_1 ~\cdots~ v_q\bigr], 1\leq q\leq n$  of sorted eigenvectors of $\Tx$. The first $q$ principal components of a vector $x=\Phi(\tilde{x})$ in the feature space are then given by $V_q^{\tr}x$. It can be shown however (see~\cite{KPCA:98}) that principal components in the
feature space can be computed in the original space with kernels using the map
 $\Pi(x) := A_q^{\tr}\mathbf{k}(x)$, where $\mathbf{k}(x) = \bigl(K(x,x_1),\ldots,K(x,x_N)\bigr)^{\tr}$.

\subsection{Model Order Reduction Map}
The method we propose consists, in essence, of collecting samples and then performing a process similar to ``simultaneous principal components analysis'' on the controllability and observability Gramian estimates in the (same) RKHS. As mentioned above, given a choice of the kernel $K$ defining a RKHS $\cH$, principal components in the feature space can be computed implicitly in the original input space using $K$. Because we will find {\em non-orthogonal} coordinates
in the feature space in which the Gramians become simultaneously diagonal, the process is not strictly speaking PCA, and the favorable properties associated with an orthonormal basis are no longer available. We will, however, continue to refer to the process of diagonalizing a covariance matrix as (K)PCA.

Turning to the controllability Gramian (the case of the observability Gramian is analogous), first note that $\widehat{W}_c$ can be viewed as the sample covariance of a collection of $N\cdot m$ vectors, scaled by $T$:
\[
\widehat{W}_c = \frac{T}{mN}\sum_{i=1}^N X(t_i)X(t_i)^{\tr} =
\frac{T}{mN}\sum_{i=1}^N\sum_{j=1}^m x^j(t_i)x^j(t_i)^{\tr}.
\nonumber\]
Thus we can form the controllability kernel matrix $K_c\in\RR^{Nm\times Nm}$ of kernel
products $(K_c)_{ij} = K(x_i, x_j)$ for $i,j=1,\ldots,Nm$ in order to carry out PCA in the feature space, where we have re-indexed the set of vectors  $\{x^{k}(t_{\ell})\}$ to use a single linear index. Similarly, we can compute the observability
kernel matrix $K_o\in\RR^{Np\times Np}$ consisting of the pairwise kernel products of the collection of data vectors described in~\eqref{eqn:obs_data}. Ordinarily, $Nm,Np\gg n$ and $K_c,K_o$ will be rank deficient.

We assume here for simplicity that the number of input excitation signals $m$ is equal to the dimension of the output $p$ so that the number of samples $N$ taken from the output and state trajectories can be the same, leading to kernel matrices $K_c$ and $K_o$ of the same size. If one adopts the set of input excitations $\{u^i(t)\}$ as above, then an alternative although more restrictive assumption can be that the number of inputs to the system is equal to the number of outputs. Then the pair $K_c,K_o$ is simultaneously diagonalized by taking the (reduced) SVD of $K_c^{1/2}K_oK_c^{1/2}$ so that $K_c^{1/2}K_oK_c^{1/2} = U\Sigma^2U^{\tr}$. Conjugation by $T = \Sigma^{1/2}U^{\tr}\sqrt{K_c^{\dagger}}$
diagonalizes $K_c$ and conjugation by $T^{-\tr}=\sqrt{\Sigma^{\dagger}}U^{\tr}K_c^{1/2}$ diagonalizes $K_o$, where $X^{\dagger}$ denotes the pseudoinverse of $X$. Finally, the order of the model is reduced by discarding small eigenvalues $\{\Sigma_{ii}\}_{i=q+1}^n$, and projecting onto the subspace associated with the first $q<n$ largest eigenvalues. This leads to the state-space reduction map $\Pi:\RR^n\to\RR^q$ given by
\begin{equation}\label{eqn:pi_map}
\Pi(x) = T_q^{\tr}\mathbf{k}_c(x), \quad x\in\RR^n
\end{equation}
where
\begin{equation}\label{eqn:emp_kmap}
\mathbf{k}_c(x) := \bigl(K(x,x^1(t_1)),\ldots,K(x,x^m(t_N))\bigr)^{\tr}.
\end{equation}


\section{Closed Dynamics of the Reduced System}\label{sec:closed_sys}
Given the nonlinear state space reduction map $\Pi:\RR^n\to\RR^q$, a remaining challenge is to
construct a corresponding (reduced) dynamical system on the reduced state space which well approximates
the input-output behavior of the original system on the original state space. Setting $x_r = \Pi(x)$ and
applying the chain rule,
\begin{equation}\label{eqn:exact_closed}
\dot{x}_r = \left.\bigl(J_{\Pi}(x)f(x,u)\bigr)\right|_{x = \Pi^{-1}(x_r)}.
\end{equation}
However we are faced with the difficulty that the map $\Pi$ is not in general injective
(even if $q=n$), and moreover one cannot guarantee that an arbitrary point in the RKHS has a non-empty
preimage under $\Phi$~\cite{mika98}. We propose an
approximation scheme to get around this difficulty: The dynamics $f$ will be approximated by an
element of an RKHS {\em defined on the reduced state space}. When $f$ is assumed to be known explicitly it can be approximated to a high degree of accuracy. An approximate, least-squares notion of ``$\Pi^{-1}$'' will be given to first or second order via a Taylor series expansion, but only where it is strictly needed -- and at the last possible moment -- so that a first or second order approximation will not be as crude as one might suppose. We will also consider, as an alternative, a direct approximation of $J_{\Pi}(\Pi^{-1}(x_r))$ which takes into account further properties of the reproducing kernel as well as the fact that the Jacobian is to be evaluated at $x = \Pi^{-1}(x_r)$ in particular. In both cases, the important ability of the map $\Pi$ to capture strong nonlinearities will not be significantly diminished.

\subsection{Representation of the dynamics in RKHS}\label{sec:f_rkhs}
The vector-valued map $f:\RR^{n}\times\RR^m\to\RR^n$ can be approximated by a composing a set of $n$ regression functions (one for each coordinate) $\hat{f}_i:\RR^{q\times m}\to\RR$ in an RKHS, with the reduction map $\Pi$. It is reasonable to expect that this approximation will be better than directly computing $f(\widehat{\Pi}^{-1}(x_r),u)$ using, for instance, a Taylor expansion notion of ``$\Pi^{-1}$'', which may ignore important nonlinearities at a stage where crude approximations must be avoided.

Let $\tilde{x}=\Pi(x)$ denote a reduced state  variable, and concatenate the input examples $\tilde{x}_j=\Pi(x_j)\in\RR^q,u_j\in\RR^m$ so that $z_j=(\tilde{x}_j,u_j)\in\RR^{q\times m}$, and $\{\left(f_i(x_j,u_j), z_j\right)\}_{j=1}^{\ell}$ is a set of input-output training pairs describing the $i$-th coordinate of the map $(\tilde{x},u)\mapsto f(x,u)$. The training examples should characterize
``typical'' behaviors of the system, and can even re-use those trajectories simulated in response to impulses for estimating the Gramians above. We will seek the function $\hat{f}_i\in\cH$ which minimizes
$$
\sum_{j=1}^{\ell}\bigl(\hat{f}_i(z_j) - f_i(x_j,u_j)\bigr)^2 + \lambda_i\|\hat{f}_i\|^2_{\cH}
$$
where $\lambda_i$ here is a regularization parameter. We have chosen the square loss, however other suitable loss functions may be used.
It can be shown~\cite{Wahba} that in this case $\hat{f}_i$ takes the form $\hat{f}_i(z) = \sum_{j=1}^{\ell}c_j^iK^f(z,z_j), i=1,\ldots,n$, where $K^f$ defines the RKHS $\cH_f$ (and is unrelated to $K$ used to estimate the Gramians). Note that although our notation takes the RKHS for each coordinate function to be the same, in general this need not be true: different kernels may be chosen for each function. Here the $\{c_j^i\}$ comprise a set of coefficients learned using the regularized least squares (RLS) algorithm. The kernel family and any hyper-parameters can be chosen by cross-validation.
For notational convenience we will further define the vector-valued empirical feature map
\begin{equation*}
\bigl(\bk^{f}(\tilde{x},u)\bigr)_i:= K^f\bigl( (\tilde{x},u), z_i \bigr)
\end{equation*}
for $i=1,\ldots,\ell$. In this notation $\hat{f}_i\bigl(\Pi(x),u\bigr) = \bc_i^{\tr}\bk^{f}(\tilde{x},u)$
where $(\bc_i)_j = c_j^i$.
%

A broad class of systems seen in the literature~\cite{scherpen_thesis} are also characterized by separable dynamics of the form $\dot{x} = f(x) + \sum_{i=1}^mg_i(x)u_i$. In this case one need
only estimate the functions $f$ and $g_i$ from examples $\{(\Pi(x_j),f(x_j))\}_j$ and $\{(\Pi(x_j),g(x_j))\}_j$.

\subsection{Approximation of the Jacobian Contribution}\label{sec:jacobian_app}
We turn to approximating the component $J_{\Pi}\bigl(\Pi^{-1}(x_r)\bigr)$ appearing in Equation~\eqref{eqn:exact_closed}.

\subsubsection{Inverse-Taylor Expansion}
A simple solution is to compute a low-order Taylor expansion of
$\Pi$ and then invert it using the Moore-Penrose pseudoinverse to obtain the approximation.
For example, consider the first order expansion $\Pi(x) \approx \Pi(a) + J_{\Pi}(a)(x-a)$. Then we
can approximate $\Pi^{-1}(x_r)$ (in the first-order, least-norm sense) as
\[
\widehat{\Pi}^{-1}(x_r) := \bigl(J_{\Pi}(a)\bigr)^{\dagger}(x_r - \Pi(a)) + a .
 \]
 We may start with $a=x_0$, but periodically update the expansion in different regions of the dynamics if desired. A good expansion point could be the estimated preimage of $x_r(t)$ returned by the algorithm proposed in~\cite{Kwok:ICML:03}.

\subsubsection{Exploiting Kernel Properties}
For certain choices of the kernel $K$ defining the Gramian feature space $\cH$, one can exploit
the fact that $K_x$ and its derivative bear a special relationship, and potentially improve the estimate for $J_{\Pi}(\Pi^{-1}(x_r))$. Perhaps the most commonly used off-the-shelf kernel families are the polynomial and Gaussian families. For any two kernels with hyperparameters $p$ and $q$ (respectively) in one of these classes, we have that $K_p = (K_q)^{p/q}$. We'll consider the polynomial kernel of degree $d$,
$K_d(x,y):=(1+\scal{x}{y})^d$ in particular; the Gaussian case can be derived using similar reasoning. For a polynomial kernel we have that
$$
\pfrac{K_d(x,y)}{x} = dK_{d-1}(x,y)y^{\tr} = d\bigl(K_d(x,y)\bigr)^{\tfrac{d-1}{d}}y^{\tr}.
$$
Recalling that $K_d(x,y)=\scal{\Phi(x)}{\Phi(y)}_{\cH}$ and $x_r=\Pi(x)=V_q^{\tr}\Phi(x)$, if $\Pi$ was invertible then we would have
$$
\left.\pfrac{K_d(x,y)}{x}\right|_{x=\Pi^{-1}(x_r)} =
d\scal{(\Phi\circ\Pi^{-1})(x_r)}{\Phi(y)}^{\tfrac{d-1}{d}}y^{\tr}.
$$
The map $\Pi$ is not injective however, and in addition the fibers of $\Phi$ may be potentially empty, so we must settle for an approximation. It is reasonable then to {\em define}
$(\Phi\circ\Pi^{-1})(x_r)$ as the solution to the convex optimization problem
\begin{equation}\label{eqn:inv-min}
\begin{aligned}
& \underset{z\in\cH}{\min} & & \nor{z}_{\cH} \\
& \text{subj. to} & & \nor{V_q^{\tr}z - x_r}_{\RR^k} = 0 .
\end{aligned}
\end{equation}
If a point $z\in\cH$ has a pre-image in $\RR^n$ this definition is consistent with composing
$\Phi$ with the formal definition $\Phi^{-1}(z) = \{x\in\RR^n~|~\Phi(x)=z\}$ and noting that in this case
$\Pi\circ\Phi^{-1} = V_q^{\tr}(\Phi\circ\Phi^{-1}) =V_q^{\tr}z$. Furthermore, a trajectory $x_r(t)$
of the closed dynamical system on the reduced statespace need not (and may not) have a counterpart in the
original statespace by virtue of the way in which ``$\Pi^{-1}$'' is used in our formulation of the reduction map and corresponding reduced dynamical system.

One will recognize that the solution $z^*$ to~\eqref{eqn:inv-min} is just the Moore-Penrose pseudoinverse $z^* = (V_q^{\tr})^{\dagger}x_r$. Inserting this solution into the feature map representation of a kernel $K$ gives the following definition for $K(\Pi^{-1}(x_r),y)$:
\begin{align*}
K(\Pi^{-1}(x_r),y) &= \scal{(\Phi\circ\Pi^{-1})(x_r)}{\Phi(y)}_{\cH} \\
 &= \scal{(V_q^{\tr})^{\dagger}x_r}{\Phi(y)}_{\cH}  = \bigl\langle x_r,V_q^{\dagger}\Phi(y)\bigr\rangle_{\RR^k} \\
 & = \bigl\langle x_r,(V_q^{\tr}V_q)^{-1}V_q^{\tr}\Phi(y)\bigr\rangle \\
 & = \bigl\langle x_r,(V_q^{\tr}V_q)^{-1}\Pi(y)\bigr\rangle \\
 & = \bigl\langle x_r,(T_q^{\tr}K_cT_q)^{-1}\Pi(y)\bigr\rangle \\
  & = \bigl\langle x_r,(T_q^{\tr}T_q\Sigma_q)^{-1}\Pi(y)\bigr\rangle .
 \end{align*}
Substituting into the derivative for a polynomial kernel $K=K_d$ gives
$$
\left.\pfrac{K_d(x,y)}{x}\right|_{x=\Pi^{-1}(x_r)} =
d\bigl\langle x_r,(T_q^{\tr}T_q\Sigma_q)^{-1}\Pi(y)\bigr\rangle^{\tfrac{d-1}{d}}y^{\tr}
$$
which immediately gives an expression for $J_{\Pi}(\Pi^{-1}(x_r))$. Note that this approximation is
global in the sense that the $q\times q$ matrix inverse  $(T_q^{\tr}T_q\Sigma_q)^{-1}$ need only  be computed once; no updating is required during simulation of the closed system.

\subsection{Reduced System Dynamics}
Given an estimate $\hat{f}\bigl(\Pi(x),u\bigr)$ of $f(x,u)$ in the RKHS $\cH_f$ and
a notion of $J_{\Pi}\bigl(\Pi^{-1}(x_r)\bigr)$ from above, we can write down a closed dynamical system on the reduced statespace. We have
\begin{align}\label{eqn:approx_xr}
\dot{x}_r &\approx \left.\bigl(J_{\Pi}(x)\hat{f}(\Pi(x),u)\bigr)\right|_{x = \Pi^{-1}(x_r)} \nonumber\\
&\approx \left.\bigl(J_{\Pi}(x)\bigr)\right|_{x = \Pi^{-1}(x_r)}\bC^{\tr}\bk^{f}(x_r,u) \nonumber\\
&\approx T_q^{\tr}J_{\mathbf{k}}\bigl(\Pi^{-1}(x_r)\bigr) \bC^{\tr}\bk^{f}(x_r,u)
\end{align}
where $\bC$ is a matrix with the vectors $\bc_i$ as its rows, and $J_{\mathbf{k}}$ is the Jacobian of the empirical feature map defined in Equation~\eqref{eqn:emp_kmap}. Here the expression 
$J_{\mathbf{k}}\bigl(\Pi^{-1}(x_r)\bigr)$ should be interpreted as notation for either of the Jacobian approximations suggested in Section~\ref{sec:jacobian_app}.


Equation~\eqref{eqn:approx_xr} is seen to give a closed nonlinear control system
expressed solely in terms of the reduced variable \mbox{$x_r\in\RR^q$}:
\begin{equation*}
\left\{
\begin{aligned}
\dot{x}_r &= T_q^{\tr}J_{\mathbf{k}}\bigl(\widehat{\Pi}^{-1}(x_r)\bigr) \bC^{\tr}\bk^{f}(x_r,u) \nonumber\\
y &= \hat{h}(x_r)
\end{aligned}
\right.
\end{equation*}
where the map $\hat{h}\circ\Pi$ modeling the output function $h:\RR^n\to\RR^p$ is estimated as described immediately below.
Although the ``true'' reduced system does not actually exist due to non-injectivity of the feature map $\Phi$, in many situations one can expect that the above system will capture the essential input-output behavior of the original system. We leave a precise analysis of the error in the approximations appearing in~\eqref{eqn:approx_xr} to future work.

\subsection{Outputs of the Reduced System}
Analogous to the case of the dynamics $f$, we are faced with two possibilities for approximating
$y = h\bigl(\Pi^{-1}(x_r)\bigr)$. We can apply the Taylor approximation $\widehat{\Pi}^{-1}$, or as in Section~\ref{sec:f_rkhs} we can estimate
a map $(\hat{h}\circ\Pi):\RR^n\to\RR^p,~x_r\mapsto y$ from the reduced state space to the output space directly, using RKHS methods. Given samples $\{\Pi(x_j),y_j\}_{j=1}^{\ell}$, each coordinate function $\bigl(\hat{h}_i\bigr)_{i=1}^p$ is given in the familiar form $\hat{h}_i(\Pi(x)) = \sum_{j=1}^{\ell} b_j^iK^h\bigl(\Pi(x),\Pi(x_j)\bigr)$,
where $K^h$ is the kernel chosen to define the RKHS, and may be different for each coordinate. It should be noted that just given the state space reduction map $\Pi$, one can immediately compare the output of the system defined by $\hat{h}(x_r)$ to the original system without defining a closed dynamics as above.
In fact with $\Pi$ and $\hat{h}$ one can design a simpler controller which takes as input the reduced state variable $x_r$, but controls the original system.

\begin{figure*}[t]
\centering
\includegraphics[height=6cm]{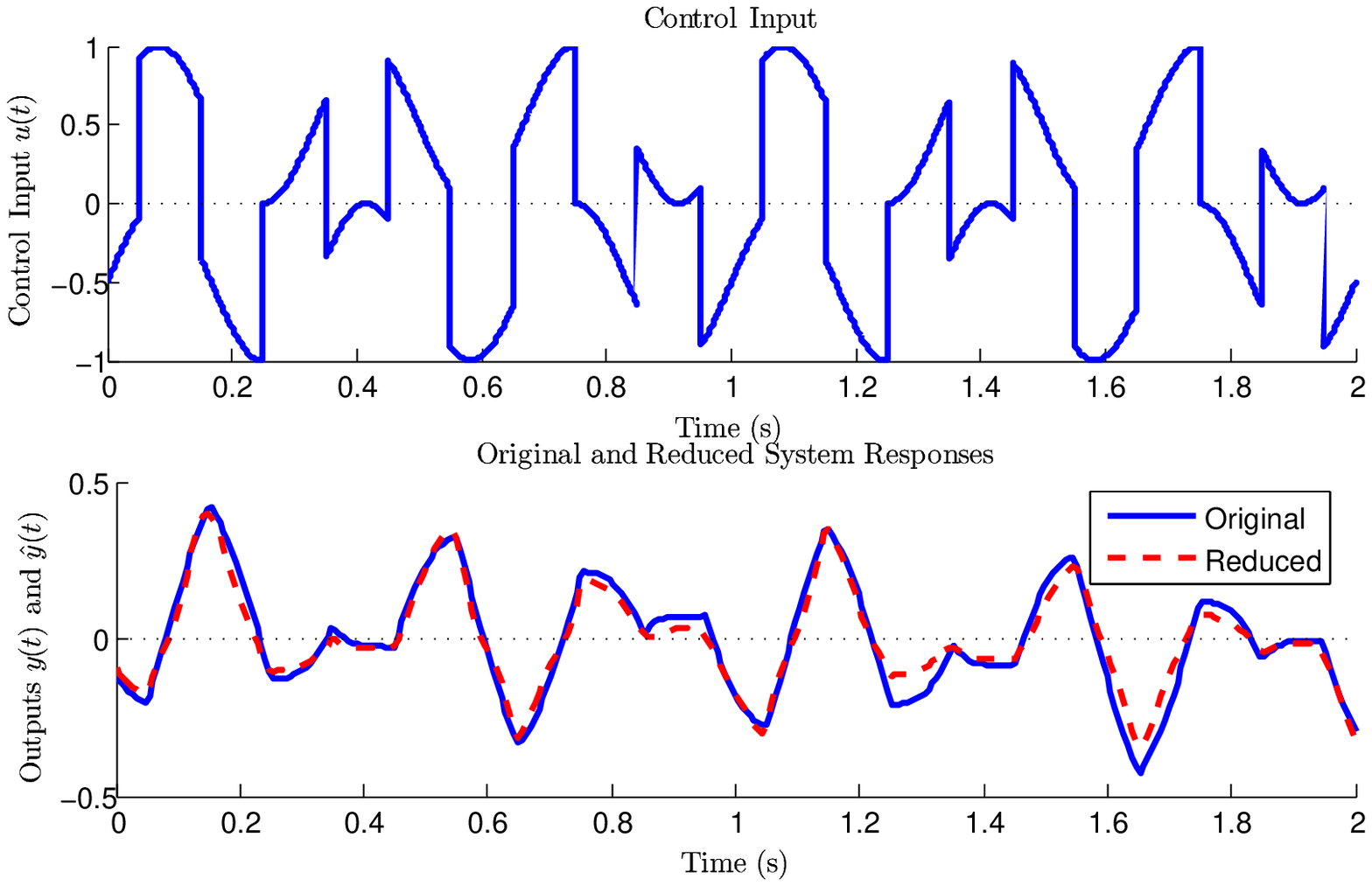}
\hskip 0.5cm
\includegraphics[height=6cm]{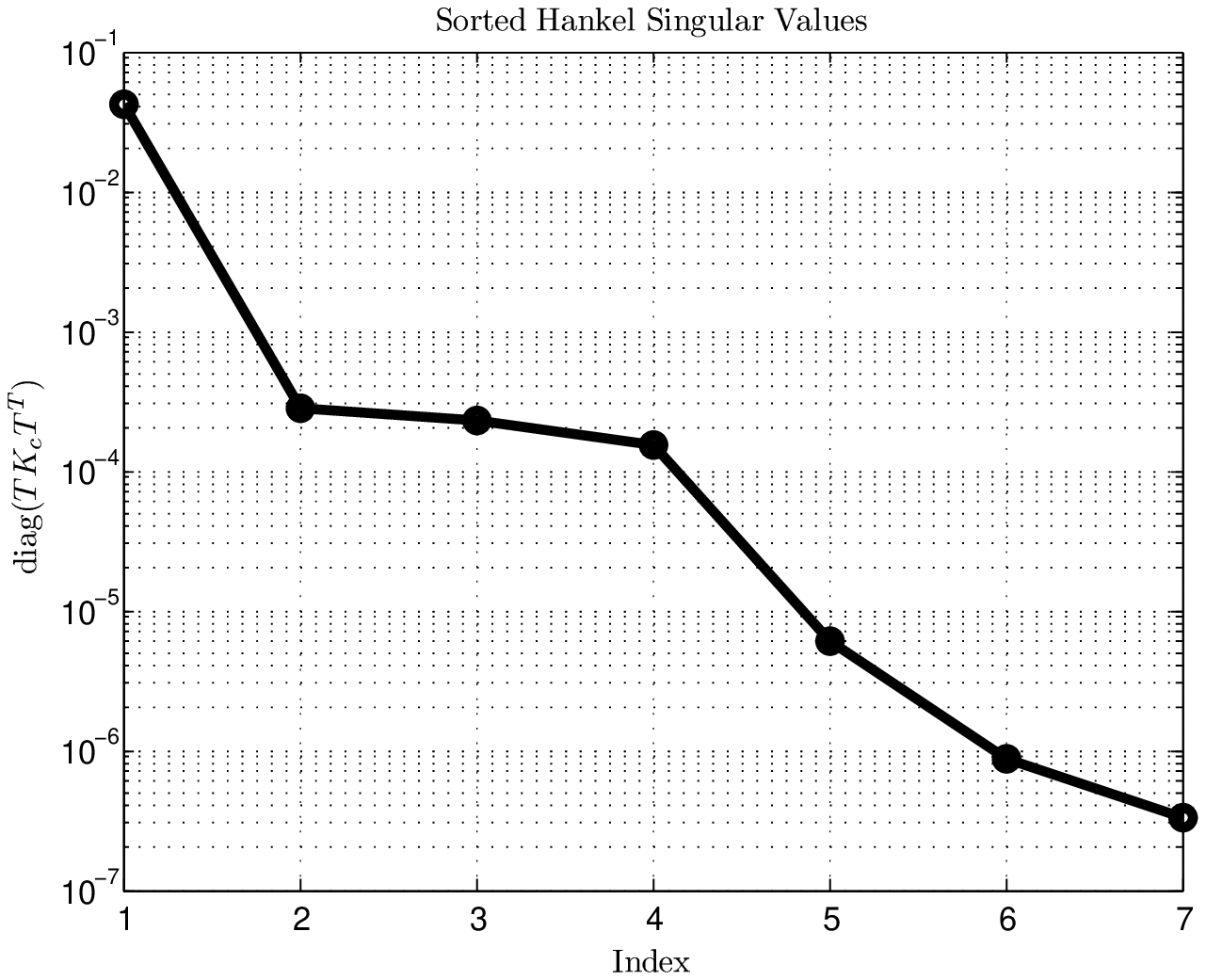}
\vskip -0.2cm
\caption{{\small\em (Left) Simulated output trajectories for the original and reduced (2-dimensional) system. (Right) Top Hankel singular values (zeros omitted).}}
\label{fig:sims}
\vskip -0.5cm
\end{figure*}

\section{Experiments}\label{sec:expts}
\begin{figure*}[t]
\centering
\includegraphics[height=5.5cm]{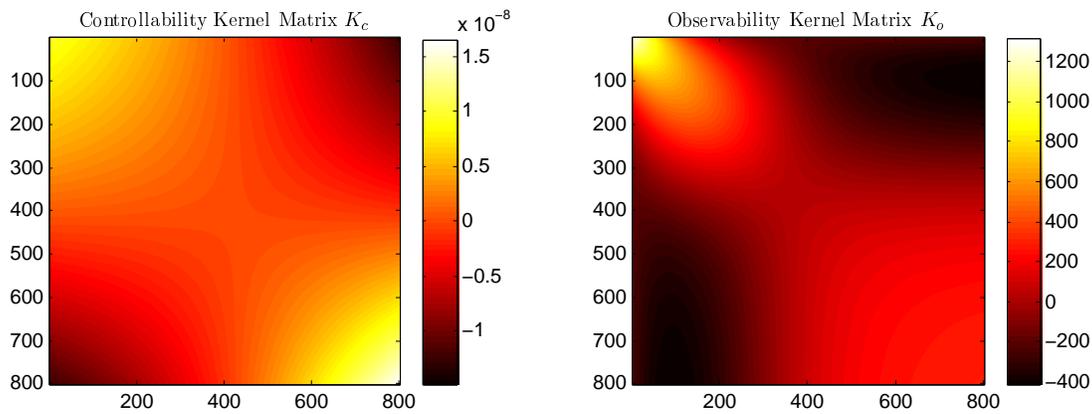}
\vskip -0.2cm
\caption{{\small\em Plots of the kernel matrices encoding information about the empirical controllability Gramian matrix (left) and
observability Gramian matrix (right).}}
\label{fig:kernmats}
\vskip -0.5cm
\end{figure*}

We demonstrate an application of our method to a 7-dimensional nonlinear system with one dimensional input and
output appearing in~\cite{nilsson} (Example 3.2, pg. 54):
\begin{gather*}
\begin{align*}
\dot{x}_1 &= -x_1^3 + u & \dot{x}_2 &= -x_2^3 -x_1^2x_2 +3x_1x_2^2 -u \\
\dot{x}_3 &= -x_3^3 + x_5 + u & \dot{x}_4 &= -x_4^3 + x_1 -x_2 +x_3 +2u \\
\dot{x}_5 &= x_1x_2x_3 -x_5^3 + u  & \dot{x}_6 &= x_5 - x_6^3 -x_5^3 +2u
\end{align*}
\\
\dot{x}_7 = -2x_6^3 +2x_5 -x_7 -x_5^3 +4u \\
y = x_1 - x_2^2 + x_3 +x_4x_3 + x_5 -2x_6 +2x_7
\end{gather*}
Impulse and initial-condition responses of the system were simulated as described above, and 800 samples
equally spaced in the time interval $[0,5s]$ were sampled to build the kernel matrices $K_c$ and $K_o$ using the third degree polynomial kernel $K(x,y)=(1 + \scal{x}{y})^3$. Recall that these kernel matrices, are the inner product counterparts to the empirical Gramians. Examples of $K_c$ and $K_o$ for this system are shown in Figure~\ref{fig:kernmats}. We imposed a small amount of
regularization when computing the balancing transformation $T$, taking the Cholesky decomposition of
$K_c + 0.001\cdot I$ instead of $K_c$. Figure~\ref{fig:sims} (right pane) shows the Hankel singular values $\Sigma=TK_cT^{\tr}$ for this problem on a log scale. One can see that perhaps the first two components ought to capture most of the system's behavior. Thus the reduction map $\Pi$ was defined by taking only the eigenvectors (scaled columns of $T$) corresponding to the largest two Hankel singular values, giving a reduced state space of dimension two.

Next, a map from the reduced variable $x_r$ to $\dot{x}$ was estimated following Section~\ref{sec:f_rkhs}.  The control input was chosen to be a 10hz square wave, and 1000 samples from the simulated system in the interval $[0,5s]$ were mapped down using $\Pi$ and then used to solve the $(n)$ RLS regression problems, one for each state variable, again using a third degree polynomial kernel. All initial conditions were set to zero. The desired outputs (dependent variable examples) used to learn $\hat{f}$ were taken to be the true $f$ evaluated at the samples from the simulated state trajectory. We also added a bias dimension of 1's to the data to account for any offset, and used a fast leave-one-out cross-validation (LOOCV) computation~\cite{RifRLS} to select the optimal regularization parameter. Two remarks are in order. The above dynamics can in fact be represented explicitly and exactly in a 3rd degree polynomial RKHS; only monomials up to degree 3 appear in the dynamics. Second, the control input is decoupled from the state. Both of these facts can be used to obtain an improved reduced model, however we did not make use of these special properties and instead applied the simplest version of the techniques described above which assume no special structure.

We followed a similar process to learn the output function $y = \hat{h}(x_r)$. Here we used a
10Hz square wave control input, zero initial conditions and 700 samples in the interval $[0,5s]$. For this function the Gaussian
kernel $K(x,y) = \exp(-\gamma\|x-y\|_2^2)$ was used to demonstrate that our method does not rely
on any particular match between the form of the dynamics and the type of kernel. The scale hyperparameter $\gamma$ was chosen to be the average distance between the training examples. We again used LOOCV
to select the RLS regularization parameter.

Finally, the closed system was simulated as described above using $x_0=0$ and a control input different from those used to learn the dynamics and output functions: $u(t) = \tfrac{1}{2}\bigl(\sin(2\pi 3t) + \text{sq}(2\pi 5t -\pi/2)\bigr)$ where $\text{sq}(\cdot)$ denotes the square wave function. This input is shown at the top of Figure~\ref{fig:sims} (left panel). The Taylor series approximation for $\Pi$ was done once about $x_0$ and was not updated further. The simulated outputs $\hat{y}(t)$ of the closed reduced system as well as the output $y(t)$ of the original system are plotted at the bottom in Figure~\ref{fig:sims} (left panel). One can see that, even for a significantly different input, the two dimensional reduced system closely captures the original system. The main source of error is seen to be
over- and under-shoot near the square wave transients. This error can be further reduced by simulating the system for different sorts of inputs (and/or frequencies) and including the collected samples in the training sets used to learn $\Pi, \hat{f}$ and $\hat{h}$. Indeed, we have had some success driving example systems with random uniform input in some cases.

\section{Conclusion}
In this paper we introduced a new model reduction method for nonlinear control systems. The method
assumes that the nonlinear system is approximately linear in a high dimensional feature space, and carries out linear balanced truncation in that space. This leads to a nonlinear reduction map, which we suggest can be combined with representations of the dynamics and output functions by elements of an RKHS to give a closed reduced order dynamical system which captures the input-output characteristics of the original system. We then demonstrated an application of our technique to a 7-dimensional system and simulated the original and reduced models for comparison, showing that the approach proposed here can yield good low-order nonlinear reductions of strongly nonlinear control systems. We believe that techniques well known to the machine learning and statistics communities can offer much to control and dynamical systems research, and many further directions remain, including reduction of unstable
systems and structure preserving systems.

\IEEEtriggeratref{9}

\end{document}